\documentclass{article}

\usepackage{amssymb}

\usepackage[english,french]{babel}
\usepackage{amsmath}
\usepackage{color}

\usepackage{bm}

\addto\captionsfrench{}

\newtheorem{theorem}{Theorem}[section]
\newtheorem{lemma}[theorem]{Lemma}
\newtheorem{e-proposition}[theorem]{Proposition}

\newtheorem{e-definition}[theorem]{Definition\rm}
\newtheorem{remark}{\it Remark\/}

\def\og{\leavevmode\raise.3ex\hbox{$\scriptscriptstyle\langle\!\langle$~}}
\def\fg{\leavevmode\raise.3ex\hbox{~$\!\scriptscriptstyle\,\rangle\!\rangle$}}

\def\eR{\mathbf{R}}
\def\eN{\mathbf{N}}

\def\dist{\operatorname{dist}}

\newcommand\dvr{\mathop{\mathrm{div}}\nolimits}
\newcommand\tder{\partial_t}

\begin{document}

\title{Extensibility of a system of transport equations in the case of an impermeable boundary}

\author{Martin Kalousek\footnote{Institute of Mathematics, University of W\"urzburg, Emil-Fischer-Str. 40, 97074 W\"urzburg, Germany, martin.kalousek@mathematik.uni-wuerzburg.de}\ , \v S\'arka Ne\v casov\'a\footnote{Institute of Mathematics, Czech Academy of Sciences, \v Zitn\'a 25, 115 67 Praha 1, Czech Republic, matus@math.cas.cz} \ and Anja Schl\"omerkemper\footnote{Institute of Mathematics, University of W\"urzburg, Emil-Fischer-Str. 40, 97074 W\"urzburg, Germany, anja.schloemerkemper@mathematik.uni-wuerzburg.de}}

\maketitle

\begin{abstract}

We show that the steady and unsteady system of transport equations with a nonhomogeneous right hand side can be extended from its domain that possesses an impermeable $C^2$-boundary to the whole space.

\vskip 0.5\baselineskip

\end{abstract}

\vspace{1cm}
For $T>0$ and a smooth bounded domain $\Omega\subset\eR^d$, $d\geq 2$, we consider the system
\begin{equation}\label{ClassForm}
		\tder h +\dvr (u\otimes h)= f \text{ in }(0,T)\times\Omega
\end{equation}
for the unknown vector function $h:(0,T)\times\Omega\to\eR^N$, the given velocity $u:(0,T)\times\Omega\to\eR^d$ and right hand side $f:(0,T)\times\Omega\to\eR^N$. We assume that the boundary $\partial\Omega$ of $\Omega$ is impermeable, i.e., denoting by $n$ the outer normal to $\partial\Omega$, the velocity $u$ satisfies
\begin{equation*}
u\cdot n=0\quad \text{ on }\partial\Omega.
\end{equation*}

The goal of this note is to investigate the extensibility of \eqref{ClassForm} from its domain $\Omega$ to the whole space. To the best of the author's knowledge there is no proof on this extension in the literature if only the normal component of the velocity is assumed to vanish on the boundary. Corresponding results for the vanishing velocity on the boundary and $N=1$ can be found in \cite[Lemma~3.3]{Feireisl-etal2001} and \cite[Lemma~6.8.]{NovStr2004}. The requirement of such extensions arises in the theory of compressible fluids. Therefore, this note can be regarded as a technical lemma needed in the proof of existence of weak solutions to the compressible Navier-Stokes equations with slip boundary conditions, see \cite[Theorem~7.69]{NovStr2004} and \cite[Theorem~3.1]{FeiNov09}. In that setting, $h$ is the mass density and $f$ represents some external force. Moreover, we expect our result to be of use in the proof of the existence of weak solutions for systems allowing for fluid-structure interaction and systems appearing in viscoelasticity, see Corollary \ref{Col:ViscoelExample}.\\

Let us start with the notation that we often use. For vectors $a\in\eR^k$, $b\in\eR^l$, $k,l\in\eN$ the outer product $a\otimes b$ denotes the $k\times l$ matrix with entries $a_ib_j$, $i=1,\ldots k$, $j=1,\ldots,l$. The inner product of vectors as well as that of matrices is denoted by $\cdot$, a centred dot. The standard notation is used for Lebesgue, Sobolev and Bochner spaces. Moreover, we denote by $W^{1,p}_0(\Omega)$ the closure of the space of smooth functions with compact support in $\Omega$ in the norm of $W^{1,p}$. Further,  $W^{1,p}_n(\Omega)$ denotes the space of functions from $W^{1,p}(\Omega)$ whose normal component of the trace vanishes on $\partial\Omega$, and $C^1_c(\Omega)$ stands for differentiable functions with compact support in $\Omega$. A generic constant is denoted by $c$.\\

In the present note we deal with a weak solution to \eqref{ClassForm} that is a function $h\in L^\infty(0,T;L^p(\Omega)^N)$, $p\in[1,\infty]$ satisfying
\begin{equation}\label{FOrigForm}
	\int_0^T\int_{\Omega}h\cdot\tder\Phi+u\otimes h\cdot\nabla\Phi+f\cdot\Phi=0\quad \text{ for all }\Phi\in C^1_c((0,T)\times\Omega)^N
\end{equation}
provided $f\in L^1(0,T;L^p(\Omega)^N)$ and $u\in L^{p'}(0,T;W^{1,{p'}}_n(\Omega)^d)$ where $p$ and $p'$ are conjugate exponents. The existence of weak solutions can be shown by adopting methods from \cite[Section IV.4]{DiPeLio89}.

The key ingredient for the proof of the main result is the following Hardy inequality that is a consequence of a general embedding for weighted Sobolev spaces, see \cite[Th\'eor\`eme 1.6]{Necas62}.
\begin{lemma}\label{lem:Hardy}
Let $\Omega\subset\eR^d$ be a bounded Lipschitz domain and $p\in(1,\infty)$. Then there is $c>0$ depending on $d,p,\Omega$ such that 
\begin{equation}\label{Hardy}
	\left\|\frac{v}{\dist(\cdot,\partial\Omega)}\right\|_{L^p(\Omega)}\leq c\|v\|_{W^{1,p}(\Omega)}\text{ for all }v\in W^{1,p}_0(\Omega).
\end{equation}
\end{lemma}
Having introduced all necessary preliminaries we can now formulate and prove the result.\\

\begin{theorem}\label{Thm:Main}
	Let $\Omega\subset\eR^d$, $d\geq 2$, be a bounded domain with a $C^2$-boundary, let $T>0$ be arbitrary and $p\in[1,\infty)$. Let $h\in L^\infty(0,T;L^p(\Omega)^N)$ be a weak solution to \eqref{ClassForm} for given $u\in L^2(0,T;W^{1,2}_n(\Omega)^d)$ and $f\in L^1(0,T; L^p(\Omega)^N)$. 
	Then, denoting by $\tilde h\in L^\infty(0,T;L^p(\eR^d)^N)$ the extension of $h$ and $\tilde f\in L^1(0,T;L^p(\eR^d)^N)$ the extension of $f$ by the zero vector in $(0,T)\times\left(\eR^d\setminus \Omega\right)$, there exists an extension $\tilde u\in L^{p'}(0,T;W^{1,p'}(\eR^d)^d)$ of $u$ such that
	\begin{equation}\label{FExtendForm}
	\int_0^T\int_{\eR^d}\tilde h\cdot\tder\Phi+\tilde u\otimes \tilde h\cdot\nabla\Phi+\tilde f\cdot\Phi=0\quad \text{ for all }\Phi\in C^1_c((0,T)\times\eR^d)^N.
	\end{equation}
\end{theorem}

\textbf{Proof:}
Let $V\subset\eR^d$ be such that $\Omega\subset\subset V$. Then there exists an extension operator $E:W^{1,p'}(\Omega)\to W^{1,p'}(\eR^d)$ which possesses the following properties, cf.\ e.g.\  \cite[Section 5.4, Theorem 1]{Evans98}:
	\begin{enumerate}
		\item $E$ is linear,
		\item $Ev=v$ a.e. in $\Omega$, and the support of $Ev$ is a subset of $V$,
		\item $\|Ev\|_{W^{1,p'}(\eR^d)}\leq c\|v\|_{W^{1,p'}(\Omega)}$ with the constant $c$ depending on $p'$, $V$ and $\Omega$ only.
	\end{enumerate}
 Let $u\in L^{p'}(0,T;W^{1,p'}_n(\Omega)^d)$ be given as in the assertion. We denote its extension for almost all $t\in(0,T)$ by  $\tilde u(t,\cdot)=Eu(t,\cdot)$. Employing the properties of $E$, we have that $\tilde u=u$ a.e. in $(0,T)\times\Omega$ and  $\tilde u\in L^{p'}(0,T;W^{1,p'}(\eR^d)^d)$. As indicated in the assertion, the extensions of $f$ and $h$  by zero are denoted by $\tilde f$ and $\tilde h$, respectively. It remains to show that $\tilde u$, $\tilde f$ and $\tilde h$ satisfy \eqref{FExtendForm}. In order to show the latter identity, we utilize a special test function in \eqref{FOrigForm} that involves the distance function $d_\Omega$ defined as
	\begin{equation*}
		d_\Omega(x) =\begin{cases}\dist(x,\partial\Omega), &x\in\Omega,\\
			-\dist(x,\partial\Omega), &x\in \eR^d\setminus\Omega.
		\end{cases}
	\end{equation*}
By \cite{Foote84,KrantzParks81} and the regularity assumptions on the boundary of $\Omega$, the distance function is $C^2$ on $V^{\sigma_0}=\{x\in\Omega: d_\Omega(x)< \sigma_0\}$ for some $\sigma_0>0$. Next, for $\sigma\in(0,\tfrac12\sigma_0)$  we consider a cut-off function $\xi^\sigma\in C^1(\eR)$ such that $0\leq \xi^\sigma\leq 1$ and
	\begin{equation*}
		\xi^\sigma(s)=\begin{cases}1& s\geq 2\sigma,\\
			0& s\leq \sigma.
		\end{cases}
	\end{equation*}
It is well-known that $\xi^\sigma$ can be chosen such that $|(\xi^{\sigma}) '|\leq c\sigma^{-1}$ for some $c>0$. Additionally, we define $\eta^\sigma(x)=\xi^\sigma(d_\Omega(x))$. This satisfies 
	\begin{equation}\label{PointCnvEta}
		\eta^\sigma\to 1 \text{ pointwisely in }\Omega\text{ as }\sigma\to 0_+
	\end{equation}
and $\eta^\sigma\in C^1_c(\Omega)$. The differentiability of $\eta^\sigma$ follows since the assumption $\partial\Omega\in C^2$ implies the differentiability of $d_\Omega$ in $V^{\sigma_0}$ and since $\eta^\sigma\equiv 1$ in $\Omega\setminus V^{2\sigma}$.	
	Fixing an arbitrary $\sigma\in(0,\tfrac12\sigma_0)$ and $\Phi\in C^1_c((0,T)\times\eR^d)^N$ we rewrite 
	\begin{align*}
		&\int_0^T \int_{\eR^d} \tilde h\cdot\tder \Phi= \int_0^T \int_{\eR^d} \tilde h\cdot \tder (\eta^\sigma\Phi) + \int_0^T\int_{\eR^d} \tilde h\cdot \tder((1-\eta^\sigma)\Phi)=I^\sigma_1+I^\sigma_2,\\
		&\int_0^T \int_{\eR^d}\tilde u\otimes \tilde h\cdot\nabla\Phi=\int_0^T \int_{\eR^d}\tilde u\otimes \tilde h\cdot\nabla(\eta^\sigma\Phi)+\int_0^T \int_{\eR^d}\tilde u\otimes \tilde h\cdot\nabla\left((1-\eta^\sigma)\Phi\right)=I^\sigma_3+I^\sigma_4,\\
		&\int_0^T\int_{\eR^d}\tilde f\cdot \Phi=\int_0^T\int_{\eR^d}\tilde f\cdot (\eta^\sigma\Phi)+\int_0^T\int_{\eR^d}\tilde f\cdot  (1-\eta^\sigma)\Phi=I^\sigma_5+I^\sigma_6.
	\end{align*}
Since $\eta^\sigma \Phi\in  C^1_c((0,T)\times\Omega)^N$, extending naturally $\eta^\sigma \Phi=0$ in $(0,T)\times\eR^d\setminus\Omega$ we have $I^\sigma_1+I^\sigma_3+I^\sigma_5=0$ due to \eqref{FOrigForm}. Hence, to conclude \eqref{FExtendForm}, we need to show that $I^\sigma_2$, $I^\sigma_4$, $I^\sigma_6$ vanish in the limit as $\sigma\to 0_+$. Applying the Lebesgue dominated convergence theorem and \eqref{PointCnvEta}, we get immediately that $I^\sigma_2$ and $I^\sigma_6$ vanish in the limit as $\sigma\to 0_+$. Expanding the derivative in $I^\sigma_4$, we obtain
\begin{equation*}
	I^\sigma_4=\int_0^T\int_{\eR^d} (1-\eta^\sigma) \tilde u\otimes\tilde h \cdot\nabla \Phi-\int_0^T\int_{\eR^d} \tilde u\otimes \tilde h \cdot \Phi\otimes\nabla\eta^\sigma=I^\sigma_{4,a}+I^\sigma_{4,b}.
\end{equation*}	
As before, $I^\sigma_{4,a} \to 0$ as $\sigma \to 0_+$ by the Lebesgue dominated convergence theorem. Furthermore, using the definition of $\eta^\sigma$ we have
\begin{align*}
	I^\sigma_{4,b}=&-\int_0^T\int_{\eR^d} (\tilde u\cdot\nabla\eta^\sigma)(\tilde h \cdot \Phi)=-\int_0^T\int_{\Omega} (u\cdot\nabla\eta^\sigma)(h \cdot \Phi).
\end{align*}
Using the notation $\Omega^\sigma=\{x\in\Omega:d_\Omega(x)\in(\sigma,2\sigma)\}$ and employing the bound $|(\xi^\sigma)'(d_\Omega)|\leq c\sigma^{-1}\leq 2cd_\Omega^{-1}$ on $\Omega^\sigma$, we estimate
\begin{align*}
	|I^\sigma_{4,b}|\leq \int_0^T\int_{\Omega^\sigma}|(\xi^\sigma)'(d_\Omega)||u\cdot\nabla d_\Omega||h||\Phi|\leq 2c\int_0^T\int_{\Omega^\sigma}\left|\frac{u\cdot\nabla d_\Omega}{d_\Omega}\right||h||\Phi|.
\end{align*}
Moreover, we consider a function $\phi\in C^1(V^{\sigma_0})$ such that $0\leq \phi\leq 1$, $\phi=1$ on $V^{\frac{1}{4}\sigma_0}$ and $\phi=0$ on $\overline{V^{\sigma_0}\setminus V^{\frac{1}{2}\sigma_0}}$. Taking into account properties of $d_\Omega$, namely $d_\Omega$ is $C^2$ on $V^{\sigma_0}$ and $\nabla d_\Omega=n$ on $\partial\Omega$, we obtain $\phi u\cdot\nabla d_\Omega\in L^{p'}(0,T;W^{1,p'}_0(V^{\frac{1}{2}\sigma_0}))$ since $u\cdot n=0$ on the boundary by assumption. Having at hand also $|\nabla d_\Omega|=1$ in $V^{\sigma_0}$, we deduce, employing the H\"older and Hardy inequalities (Lemma \ref{lem:Hardy}),
\begin{align*}
	|I^\sigma_{4,b}|\leq& c\|h\|_{L^p(0,T;L^p(\Omega^\sigma))}\left\|\frac{u\cdot\nabla d_\Omega}{d_\Omega}\right\|_{L^{p'}(0,T;L^{p'}(V^{\sigma_0}))}\\
	\leq& c\|h\|_{L^p(0,T;L^p(\Omega^\sigma))}\left(\left\|\frac{\phi u\cdot\nabla d_\Omega}{d_\Omega}\right\|_{L^{p'}(0,T;L^{p'}(V^{\frac{1}{2}\sigma_0}))}+\left\|\frac{(1-\phi) u\cdot\nabla d_\Omega}{d_\Omega}\right\|_{L^{p'}(0,T;L^{p'}(V^{\sigma_0}\setminus V^{\frac{1}{4}\sigma_0}))}\right)\\
	\leq& c\|h\|_{L^p(0,T;L^p(\Omega^\sigma))}\left(\left\|\phi u\cdot\nabla d_\Omega\right\|_{L^{p'}(0,T;W^{1,p'}(V^{\frac{1}{2}\sigma_0}))}+\frac{4}{\sigma_0}\|u\|_{L^{p'}(0,T;L^{p'}(V^{\sigma_0}\setminus V^{\frac{1}{4}\sigma_0}))}\right),
\end{align*}
where the constant $c$ depends on $d$, $p$, $\sigma_0$, $\Omega$ and $\|\Phi\|_{L^\infty((0,T)\times\eR^d)}$. Therefore,  $|I^\sigma_{4,b}|\to 0$ since $|\Omega^\sigma|\to 0$ and $\|h\|_{L^p(0,T;L^p(\Omega^\sigma))}\to 0$ accordingly as $\sigma\to 0_+$. Hence the theorem is proved. \hfill $\Box$\\

\begin{remark}
The extension for the steady variant of \eqref{ClassForm} can be performed exactly in the same way as in the proof of Theorem \ref{Thm:Main}. Indeed, the only difference is the lack of the term with the time derivative of a test function; the other terms are handled as above.\\
\end{remark}

\begin{remark}\label{Col:ViscoelExample} The right hand side of \eqref{ClassForm} can have a more general form. For instance, $f(t,x,h)=f_1(t,x)+A(t,x)h$ where $A$ stands for an $N\times N$ matrix valued function. \ Such a right hand side is for instance observed in models from viscoelasticity, see, e.g., \cite{LiuWalkington}. Indeed, for the matrix valued deformation gradient $F$ one considers 
\begin{equation*}
	\tder F^m+\dvr(u\otimes F^m)=(\nabla uF)^m, \quad m=1,\ldots, d,
\end{equation*}
in the compressible as well as in the incompressible case provided the functions $u$ and $F$ enjoy  an appropriate regularity. Here the superscript $m$ stands for the $m$-th column of the corresponding matrix. The extension procedure from Theorem \ref{Thm:Main} can be easily adopted for such a system.
\end{remark}

\section*{Acknowledgements}
The authors would like to thank E.~Feireisl and A.~Novotn\'y for helpful discussions about the extension of a transport equation from its domain to the whole space. A.S.\ and M.K. were partially supported by DFG grant SCHL 1706/4-1. The research of  \v{S}.N.\  was supported by GA\v{C}R project P201-16-032308
 and RVO 67985840. Part of this work was carried out when \v S.N.\ was engaged as a Giovanni Prodi Chair Professor at the University of W\"urzburg.


\begin{thebibliography}{99}
\bibitem{DiPeLio89}
DiPerna~R.J. and Lions~P.-L., Ordinary differential equations, transport theory and {S}obolev spaces,  \emph{Invent. Math.} \textbf{98} (1989), no.~3, 511--547.
\bibitem{Evans98}
Evans~L.C., \emph{Partial differential equations}, Graduate Studies in
  Mathematics, vol.~19, American Mathematical Society, Providence, RI, 1998.

\bibitem{Feireisl-etal2001}
Feireisl~E., Novotn\'{y}~A. and Petzeltov\'{a}~H., On the existence of globally defined weak solutions to the {N}avier-{S}tokes equations,  \emph{J. Math. Fluid Mech.} \textbf{3} (2001), no.~4, 358--392.


\bibitem{FeiNov09}
Feireisl~E. and Novotn\'{y}~A., \emph{Singular limits in
  thermodynamics of viscous fluids}, Advances in Mathematical Fluid Mechanics,
  Birkh\"{a}user Verlag, Basel, 2009.

\bibitem{Foote84}
Foote~R.L., Regularity of the distance function,  \emph{Proc. Amer. Math.
  Soc.} \textbf{92} (1984), no.~1, 153--155.
\bibitem{KrantzParks81}
Krantz~S.G. and Parks~H.R., Distance to {$C\sp{k}$} hypersurfaces,  \emph{J. Differential Equations} \textbf{40} (1981), no.~1, 116--120.
\bibitem{LiuWalkington}
Liu~C. and Walkington~N.J., An {E}ulerian Description of Fluids Containing Visco-Elastic Particles, \emph{Arch. Rational Mech. Anal.} \textbf{159} (2001), 229--252.
\bibitem{Necas62}
Ne\v{c}as~J., Sur une m\'ethode pour r\'esoudre les \'equations aux d\'eriv\'ees partielles du type elliptique, voisine de la variationnelle, \emph{Ann. Scuola Norm. Sup. Pisa} \textbf{16} (1962), no.~3, 305--326.
\bibitem{NovStr2004}
Novotn\'{y}~A. and Stra\v{s}kraba~I., \emph{Introduction to the mathematical
  theory of compressible flow}, Oxford Lecture Series in Mathematics and its
  Applications, vol.~27, Oxford University Press, Oxford, 2004.
\end{thebibliography}
\end{document}